%&amstex          
\input amstex\documentstyle{amsppt}  
\pagewidth{12.5cm}\pageheight{19cm}\magnification\magstep1
\topmatter
\title{Constructible representations and Catalan numbers}\endtitle
\author George Lusztig and Eric Sommers\endauthor
\address{Department of Mathematics, M.I.T., Cambridge, MA 02139;
Department of Mathematics, University of Massachusetts, Amherst,
MA 01003}\endaddress
\dedicatory{Dedicated to the memory of Gary Seitz}\enddedicatory

\thanks{G.L. supported by NSF grant DMS-2153741}\endthanks

\endtopmatter   
\document

\define\si{\sim}

\define\sqc{\sqcup}

\define\lb{\linebreak}

\define\op{\oplus}
   
\define\part{\partial}
\define\emp{\emptyset}

\define\m{\mapsto}
\define\do{\dots}

\define\sub{\subset}    

\define\T{\times}
\define\ti{\tilde}
\define\nl{\newline}

\define\fra{\frac}

\define\e{\epsilon}

\define\s{\sigma}

\define\x{\xi}

\define\Th{\Theta}

\define\CC{\bold C}

\define\NN{\bold N}

\define\QQ{\bold Q}

\define\ZZ{\bold Z}

\define\cb{\Cal B}
\define\cc{\Cal C}

\define\ce{\Cal E}
\define\cf{\Cal F}
\define\cg{\Cal G}

\define\cl{\Cal L}

\define\ct{\Cal T}

\define\cz{\Cal Z}

\define\ta{\ti a}
\define\tb{\ti b}

\define\tE{\ti E}

\define\tL{\ti L}

\head 0. Introduction\endhead
\subhead 0.1\endsubhead
The sequence of Catalan numbers is the sequence
$Cat_n,(n=1,2,3,\do)$ where $Cat_n=\fra{(2n)!}{n!(n+1)!}$.
According to \cite{La}, Catalan numbers first appeared in the
work of Ming Antu (1692-1763). They were rediscovered by
Euler (1707-1783). See also \cite{St}.

In this paper we give a new way in which Catalan numbers appear
in connection with Lie theory.

\subhead 0.2\endsubhead
Let $G$ be a connected reductive algebraic group of adjoint type
over $\CC$ whose Weyl
group $W$ is assumed to be irreducible. Let $\hat W$ be the set of
(isomorphism classes of) irreducible representations (over $\QQ$) of $W$.

In \cite{L79}, a partition of $\hat W$ into subsets called {\it families}
was defined and in \cite{L82} a class of not necessarily
irreducible representations (later called
{\it constructible representations}, see \cite{L03}) of $W$ with all components in
a family $c$ (which we now fix) was defined by an inductive procedure.
Let $Con(c)$ be the set of constructible representations (up to
isomorphism) attached to $c$. In \cite{L82} it was
conjectured that the representations in $Con(c)$ are precisely the
representations associated in \cite{KL} to the various left cells of
$W$ contained in the two-sided cell of $W$ defined by $c$; this 
conjecture was proved in \cite{L86}. 
It is known that $|c|=1$ if $W$ is of type $A$,
$|c|=\binom{D+1}{D/2}$ (with $D\in2\NN$) if $W$ is of type
$B,C$ or $D$, and $|c|$ is one of $1,2,3,4,5,11,17$ if $W$ is
of exceptional type.

\subhead 0.3\endsubhead
We would like to find an explicit formula for $|Con(c)|$.

If $|c|$ is one of $1,2,3,4,5,11,17$ then $|Con(c)|$ is
$1,1,2,2,3,5,7$ respectively.

{\it In the remainder of this paper we
assume that }

(a) $|c|=\binom{D+1}{D/2} \text{ with } D=2d\in2\NN$.
\nl
In \S1 we prove the following result.
\proclaim{Theorem 0.4} We have $|Con(c)|=Cat_{d+1}$.
\endproclaim
It is known (see \cite{L22,2.13}) that if $W$ is of type $D$ then
$|Con(c)|=|Con(c')|$ for some family $c'$ in a Weyl group of type $B$ or
$C$. We will therefore assume in the rest of the paper
that $W$ is of type $B$ or $C$.

\subhead 0.5\endsubhead
According to \cite{HM, Cor.4}, we have
$$Cat_n=\sum_{p=1}^nN(n,p)\tag a$$
where
$$N(n,p)=(1/n)\binom{n}{p}\binom{n}{p-1}$$
are the Narayana numbers.

We denote by $F$ the field with two elements.

In \cite{L87} a bijection between $Con(c)$ and a certain
collection $X_c$ of subgroups of $F^d$ is described. 
For each $p$, $1\le p\le d+1$ let $X_{c,p}$ be the set of
subgroups of cardinal $2^{p-1}$ in $X_c$.
Th following refinement of Theorem 0.4 is proved in \S2.
\proclaim{Theorem 0.6} We have $|X_{c,p}|=N_{d+1,p}$.
\endproclaim

\subhead 0.7\endsubhead
In \S3 we state a conjecture according to which Catalan numbers appear
in connection with the study of Springer fibres for $G$.

\subhead 0.8\endsubhead
For any $i\le j$ in $\ZZ$ we set $[i,j]=\{h\in\ZZ;i\le h\le j\}$.

\head 1. Proof of Theorem 0.4\endhead
\subhead 1.1\endsubhead
Let $D\in2\NN$.
Let $V_D$ be an $F$-vector space with a nondegenerate
symplectic form $<,>:V_D\T V_D@>>>F$ and with a given
subset $\{e_1,e_2,e_3,\do,e_D\}$ such that
$<e_i,e_j>=1$ if $i-j=\pm1$  and $<e_i,e_j>=0$ otherwise.

Assuming that $D\ge2$ and $i\in[1,D]$ we define a linear (injective)
map $T_i:V_{D-2}@>>>V_D$ by 

$e_a\m e_a$  if $a<i-1$,

$e_{i-1}\m e_{i-1}+e_i+e_{i+1}$,

$e_a\m e_{a+2}$ if $a\ge i$.

(We regard $V_{D-2}$ as a subspace of $V_D$ in an obvious way.)

Let $\cf(V_D)$ be the family of isotropic subspaces associated
in \cite{L20,1.17} to $V_D$ and its basis $\{e_1,e_2,\do,e_D\}$.
(The characteristic functions of these subspaces form a basis
of the $\CC$-vector space of functions $V_D@>>>\CC$.)
We have a partition $\cf(V_D)=\sqc_{k\ge0}\cf^k(V_D)$.
We will only give here the definition of $\cf^0(V_D)$
and $\cf^1(V_D)$. The definition is by induction on $D$. When $D=0$,
$\cf^0(V_D)$ consists of $0$ and $\cf^1(V_D)$ is empty. Assume now
that $D\ge2$. A subspace
$E$ of $V_D$ is said to be in $\cf^0(V_D)$ if either $E=0$
or if there exists $i\in[1,D]$ and $E'\in\cf^0(V_{D-2})$ such that
$E=T_i(E')+Fe_i$.
A subspace $E$ of $V_D$ is said to be in $\cf^1(V_D)$ if either
$E=F(e_1+e_2+\do+e_D)$ or if there exists $i\in[1,D]$ and
$E'\in\cf^1(V_{D-2})$ such that $E=T_i(E')+Fe_i$.

For example if $D=2$, $\cf^0(V_D)$ consists of $0,Fe_1,Fe_2$ and
$\cf^1(V_D)$ consists of $F(e_1+e_2)$. If $D=4$, $\cf^0(V_D)$ consists
of
$$\align&0,Fe_1,Fe_2,Fe_3,Fe_4,Fe_1+Fe_3,Fe_1+Fe_4,Fe_2+Fe_4,\\&
F(e_1+e_2+e_3)+F(e_2),F(e_2+e_3+e_4)+F(e_3)\endalign$$
and  $\cf^1(V_D)$ consists of
$$\align&F(e_1+e_2+e_3+e_4), F(e_1+e_2+e_3+e_4)+Fe_2,\\&
F(e_1+e_2+e_3+e_4)+Fe_3,F(e_1+e_2)+Fe_4,Fe_1+F(e_3+e_4).\endalign$$
We have $$\cf^0(V_D)=\cf_{D/2}^0(V_D)\sqc\cf_{<D/2}^0(V_D)$$
where $$\cf_{D/2}^0(V_D)=\{E\in\cf^0(V_D);\dim(E)=D/2\},$$
$$\cf_{<D/2}^0(V_D)=\{E\in\cf^0(V_D);\dim(E)<D/2\}.$$

\subhead 1.2\endsubhead
Let $\cg^0_D$ (resp. $\cg_D^1$) be the set of lines in $V_D$ of the
form $F(e_a+e_{a+1}+\do+e_b)$ where $a\le b$ in $[1,D]$
satisfy $b-a=1\mod2$ (resp. $b-a=0\mod2$). Let
$\cg_D=\cg^0_D\sqc\cg^1_D$. 
For $E\in\cf(V_D)$ let $B_E=\{L\in\cg_D;L\sub E\}$.
According to \cite{L22, 1.2(e),(f),(g)},
if $E\in\cf(V_D)$ then $E=\op_{L\in B_E}L$;
moreover we have $E\in\cf^0(V_D)$ if and only if
$B_E\sub\cg^1_D$; we have $E\in\cf^1(V_D)$ if and only if
$B_E$ contains a unique line $L_E$ in $\cg^0_D$.

It follows that if $E\in\cf^1(V_D)$ we can write
$E=E_0+L_E$ where $E_0=\op_{L\in B_E;L\ne L_E}L$.

We show:

(a) $E_0\in\cf^0(V_D)$.
\nl
We argue by induction on $D$.
If $D=0$ then $\cf^1_D=\emp$ and there is nothing to prove.
Assume now that $D\ge2$. If $E=F(e_1+e_2+\do+e_D)$, then
$E_0=0$ and (a) is obvious. If $E$ is not of this form then there
exists $i\in[1,D]$ and $E'\in\cf^1_{D-2}$ such that
$E=T_i(E')+Fe_i$. By the induction hypothesis we have
$E'=E'_0\op L_{E'}$ where $E'_0\in\cf^0_{D-2}$. We have
$E=T_i(E'_0)+Fe_i+T_i(L_{E'})=\tE_0+\tL$ where 
$\tE_0=T_i(E'_0)+Fe_i\in\cf^0(V_D)$ and $\tL=T_i(L_{E'})\in\cg^0_D$
(from the definition of $T_i$). Since $\tL\sub E$ we must have
$\tL=L_E$. We have $B_E=B_{\tE_0}\cup\{L_E\}$ (the union is disjoint
since $B_{\tE_0}\sub\cg^1_D,L_E\in\cg^0_D$. Thus
$B_{\tE_0}=B_E-\{L_E\}$.
Since $\tE_0=\sum_{L\in B_{\tE_0}}=\sum_{L\in B_E-\{L_E\}}L=E_0$
we see that $E_0=\tE_0\in\cf^0(V_D)$. This proves (a).

Note that in (a) (which is a direct sum)
we have $\dim(E)\le D/2$, $\dim(L_E)=1$ hence
$\dim(E_0)<D/2$. Thus we can define a map
$\Xi_D:\cf^1(V_D)@>>>\cf^0_{<D/2}(V_D)$ by $E\m E_0$ (notation
of (a)).

We show:

(b) For any $E_0\in\cf^0_{<D/2}(V_D)$ there exists
$E\in\cf^1(V_D)$ such that $\Xi_D(E)=E_0$.
\nl
We argue by induction on $D$.
If $D=0$ then $\cf^0_{<D/2}(V_D)$ is empty 
and there is nothing to prove.
Assume now that $D\ge2$. If $E_0=0$ then $E=F(e_1+e_2+\do+e_D)$
is as required. Now assume that $E_0\ne0$. Then there exists
$i\in[1,D]$ and $E'_0\in\cf^0(V_{D-2})$ such that
$E_0=T_i(E'_0)+Fe_i$.
Since this sum is necessarily a direct sum we have
$\dim(E'_0)=\dim(T_i(E'_0))=\dim(E_0)-1<(D/2)-1=(D-2)/2$
so that $E'_0\in\cf^0_{<(D-2)/2})V_{D-2}$.
By the induction hypothesis there exists
$L\in\cg^0_{D-2}$ such that $E'_0+L\in\cf^1(V_{D-2})$.
Let $E=T_i(E'_0+L)+Fe_i$. We have $E\in\cf^1(V_D)$ and $E=E_0+T_i(L)$.
Note that $T_i(L)\in\cg^0_D$ and is contained in $E$ hence it
is equal to $L_E$.
It follows that $E_0=\Xi_D(E)$. This proves (b).

We show:

(c) Assume that $E,E'$ in $\cf^1(V_D)$ satisfy $\Xi(E)=\Xi(E')$.
Then $E=E'$.
\nl
We have $E=E_0\op L,E'=E_0\op L'$ where $E_0\in\cf^0(V_D)$
and $L=F(e_a+e_{a+1}+\do+e_b)$,
$L'=F(e_{a'}+e_{a'+1}+\do+e_{b'})$, where $a<b$ in $[1,D]$
$a'<b'$ satisfy $b-a=1\mod2$, $b'-a'=1\mod2$. (In fact, from
\cite{L20, 1.3(e), see $(P_2)$} we have that
$a=1\mod2,b=0\mod2,a'=1\mod2,b'=0\mod2$.)
Assume first that $a<a'$ so that $a\le a'-2$.
From \cite{L20, 1.3(e), see $(P_2)$} we see that
there exist $1\le c\le c'\le D$ such that $c\le a\le c'$ and such
that the line $\cl=F(e_c+e_{c+1}+\do+e_{c'})$ is contained in
$E_0$ hence also in $\cg^1_D$. But then the pair of distinct lines
$\cl,L$ would violate \cite{L20, 1.3(e), see $(P_0)$}.
We see that we must have $a\ge a'$. Similarly we have $a'\ge a$
hence $a'=a$.

Assume next that that $b<b'$ so that $b+2\le b'$.
From \cite{L20, 1.3(e), see $(P_2)$} we see that
there exist $1\le c\le c'\le D$ such that $c\le b'\le c'$ and such
that the line $\cl=F(e_c+e_{c+1}+\do+e_{c'})$ is contained in
$E_0$ hence also in $\cg^1_D$. But then the pair of distinct lines
$\cl,L'$ would violate \cite{L20, 1.3(e), see $(P_0)$}.
We see that we must have $b\ge b'$. Similarly we have $b'\ge b$
hence $b'=b$.

We see that $L=L'$ hence $E=E'$. This proves (c).

\subhead 1.3\endsubhead
From (a),(b),(c) we see that

$|\cf^0_{<D/2}(V_D)|=|\cf^1(V_D)|$
\nl
hence $|\cf^0(V_D)|-|\cf^0_{D/2}(V_D)|=|\cf^1(V_D)|$ that is,
$$|\cf^0_{D/2}(V_D)|=|\cf^0(V_D)|-|\cf^1(V_D)|.$$
According to \cite{L20, 1.27} we have
$$|\cf^0(V_D)|=\binom{D+1}{D/2}, |\cf^1(V_D)|=\binom{D+1}{(D-2)/2}.$$
It follows that 
$$|\cf^0_{D/2}(V_D)|=\binom{D+1}{D/2}-\binom{D+1}{(D-2)/2}
=\fra{(2d+2)!}{(d+1)!(d+2)!}=C_{d+1}$$
where $D=2d$.

\subhead 1.4\endsubhead
In \cite{L81} the set $c$ is identified with a subset
of $V_D$. Now any object in $Con(c)$ is multiplicity free hence
may be identified with a subset of $c$ hence with a subset
of $V_D$. This subset is a Lagrangian subspace of $V_D$.
Thus $Con(c)$ is identified with a subset of the set of
Lagrangian subspaces of $V_D$.
This subset is the same as $\cf^0_{D/2}(V_D)$ (see
\cite{L19, 2.8(iii)}).
We see that $|Con(c)|=C_{d+1}$ and Theorem 0.4 is proved.

\head 2. Proof of Theorem 0.6\endhead
\subhead 2.1\endsubhead
We preserve the notation of $V_D$. We have $V_D=V_D^0\op V_D^1$ where
$V_D^0$ has basis $\{e_2,e_4,\do,e_D\}$ and $V_D^1$ has basis
$\{e_1,e_3,\do,e_{D-1}\}$.
Assuming that $D\ge2$ we define for any $i\in[1,D]$
a linear map $\ct_i:V_{D-2}^1@>>>V_D^1$ by

$e_k\m e_k$ if $k\le i-2$,

$e_k\m e_{k+2}$ if $k\ge i$,

$e_{i-1}\m\{e_{i-1},e_{i+1}\}$ if $i$ even.

Following \cite{L19, 2.3} we define a collection $\cc(V^1_D)$
of subspaces of $V_D^1$ by induction on $D$. If $D=0$, $\cc(V^1_D)$
consists of $\{0\}$.
Assume now that $D\ge2$. A subspace $\ce$ of $V_D^1$ is said
to be in $\cc(V^1_D)$ if either $\ce=\{0\}$ or if there exists
$i\in[1,D]$ and $\ce'\in \cc(V^1_{D-2})$ such that

$\ce=\ct_i(\ce')+Fe_i$ (if $i$ is odd)

$\ce=\ct_i(\ce')$ (if $i$ is even).

For example,  $\cc(V^1_2)$ consists of $2$ subspaces: $0,Fe_1$;
 $\cc(V^1_4)$ consists of $5$ subspaces:

$0, Fe_1,Fe_3,F(e_1+e_3),Fe_1+Fe_3$;

 $\cc(V^1_6)$ consists of $14$ subspaces:

$0, Fe_1,Fe_3,Fe_5,F(e_1+e_3),F(e_3+e_5),F(e_1+e_3+e_5)$,

$Fe_1+Fe_3,Fe_1+Fe_5,Fe_3+Fe_5,F(e_1+e_3)+Fe_5,Fe_1+F(e_3+e_5),
F(e_1+e_3+e_5)+Fe_3,Fe_1+Fe_3+Fe_5$.

\subhead 2.2\endsubhead
If $\ce\in \cc(V^1_D)$ we set $\ce^!=\{x\in V_D^0;<x,\ce>=0\}$. 
The following result appears in \cite{L19, 2.4}.

(a) $\ce\m\ce\op\ce^!$ defines a bijection 
$\cc(V^1_D)@>\si>>\cf^0_{D/2}(V_D)$. The inverse bijection is
given by $E\m E\cap V^1_D$.

\subhead 2.3\endsubhead
Let $\cz^*_D$ be the set of all elements of $V^1_D$ of the form

$e_{a,b}=e_a+e_{a+2}+e_{a+4}+\do+e_b$

for various numbers $a\le b$ in $\{1,3,\do,D-1\}$.

For any $s\ge0$ let
$\cz^s_D$ be the set of all finite unordered sequences
$$e_{a_1,b_1},e_{a_2,b_2},\do,e_{a_s,b_s}$$ in $\cz^*_D$ such that
for any $n\ne m$ in $\{1,2,\do,s\}$ we have either

$a_n\le b_n<a_m\le b_m$ or $a_m\le b_m<a_n\le b_n$,

or $a_n<a_m\le b_m<b_n$ or $a_m<a_n\le b_n<b_m$.

Let $\cz_D=\cup_{s\ge0}\cz^s_D$ (a disjoint union).

For example, $\cz_2$ consists of $2$ sequences: $\emp, \{e_1\}$;

$\cz_4$ consists of $5$ sequences:
$\emp, \{e_1\},\{e_3\},\{e_1+e_3\},\{e_1,e_3\}$;

$\cz_6$ consists of $14$ sequences:

$\emp,\{e_1\},\{e_3\},\{e_5\},\{e_1+e_3\},\{e_3+e_5\},
\{e_1+e_3+e_5\}$,

$\{e_1,e_3\},\{e_1,e_5\},\{e_3,e_5\},\{e_1+e_3,e_5\}\{e_1,e_3+e_5\},
\{e_1+e_3+e_5,e_3\},\{e_1,e_3,e_5\}$.

We have the following result.
\proclaim{Theorem 2.4} The assignment
$$\Th_D:(e_{a_1,b_1},e_{a_2,b_2},\do,e_{a_s,b_s})\m
Fe_{a_1,b_1}+Fe_{a_2,b_2}+\do+Fe_{a_s,b_s}$$
defines a bijection $\cz_D@>\si>>\cc(V^1_D)$.
\endproclaim
When $D\le6$ this follows from 2.1, 2.3.
Note that the Theorem gives an order preserving bijection between
the set of non crossing partitions (see \cite{St}) and
$\cc(V^1_D)$ (with the order given by inclusion).

\subhead 2.5\endsubhead
Assuming that $D\ge2$ we define for any $i\in[1,D]$
a map $\s_i:\cz^*_{D-2}@>>>\cz^*_D$ by

$e_{a,b}\m e_{a+2,b+2}$ if $i\le a$,

$e_{a,b}\m e_{a,b+2}$ if $a<i\le b+1$, 

$e_{a,b}\m e_{a,b}$ if $i>b+1$. 

Note that

$\s_i(e_{a,b})=\ct_i(e_{a,b})$ if $i$ is even,

$\s_i(e_{a,b})=\ct_i(e_{a,b})$ if $i$ is even and $i\le a$ or $i>b$,

$\s_i(e_{a,b})=\ct_i(e_{a,b})+e_i$ if $i$ is odd and $a<i\le b$.

\subhead 2.6\endsubhead
Assume that $D\ge2$ and $i\in[1,D]$. Let $e_{a,b},e_{a',b'}$ be in
$\cz^*_{D-2}$ and let $e_{\ta,\tb}=\s_i(e_{a,b})$,
$e_{\ta',\tb'}=\s_i(e_{a',b'})$. We show:

(i) If $b<a'$ then $\tb<\ta'$.

(ii) If $a<a'$ and $b'<b$ then $\ta<\ta'$ and $\tb'<\tb$. 

(iii) If $i$ is odd and $\ta\le i\le\tb$ then $\ta<i<\tb$.
\nl
In the setup of (i) assume that $\ta'\le\tb$. Then we have
$a'\le b$ or $a'+2\le b$ or $a'+2\le b+2$ or $a'\le b+2$. The first
3 cases are clearly impossible; in the 4th case we have $b+2=a'$
(since $b+2\le a'\le b+2$), $b'+1<i$ and $b+1\ge i$, so that
$b>b'\ge a'$, a contradiction.

In the setup of (ii) assume that $\ta\ge\ta'$. Then we have
$a\ge a'$ or $a+2\ge a'+2$ or $a\ge a'+2$ or $a+2\ge a'$.
The first 3 cases are clearly impossible, in the 4th case we have
$a+2=a'$ (since $a+2\le a'\le a+2$), $a'<i$ and $a\ge i$, so that
$a>a'$, a contradiction. Thus, $\ta<\ta'$.

Again, in the setup of (ii) assume that $\tb'\ge\tb$. Then we  have
$b'\ge b$ or $b'+2\ge b+2$ or $b'\ge b+2$ or $b'+2\ge b$. The first 3
cases are clearly impossible. In the 4th case we have
$b'+2=b$ (since $b\ge b'+2\ge b$), $b+1<i$ and $b'+1\ge i$ so that
$b'>b$, a contradiction. Thus, $\tb'<\tb'$.

In the setup of (iii) assume that $\ta=i$. We have $\ta=a$ or $\ta=a+2$.
If $\ta=a$ we have $a=i$ and $b<i$ hence $b<\tb$ so that $\tb=b+2$;
this implies $i\le b$, a contradiction. 
If $\ta=a+2$ we have $a+2=i$, $i\le a$, a contradiction. Thus $\ta<i$.

In the setup of (iii) assume that $\tb=i$.
We have $\ta=b$ or $\tb=b+2$.
If $\tb=b$ we have $b=i$ and $b<i$, a contradiction.
If $\tb=b+2$ we have $b+2=i$ and either $a\ge i$ or $a<i\le b$.
In the first case we have $a\ge b+2>b$, a contradiction; in the
second case we have $b+2\le b$,  a contradiction. Thus, $i<\tb$.

\subhead 2.7\endsubhead
From 2.6(i)-(iii) we see that when $D\ge2$ and $i\in[1,D]$,
there is a well defined map $\Sigma_i:\cz_{D-2}@>>>\cz_D$ given by
$$(e_{a_1,b_1},e_{a_2,b_2},\do,e_{a_s,b_s})\m
(\s_i(e_{a_1,b_1}),\s_i(e_{a_2,b_2}),\do,\s_i(e_{a_s,b_s}),e_i)$$
if $i$ is odd,
$$(e_{a_1,b_1},e_{a_2,b_2},\do,e_{a_s,b_s})\m
(\s_i(e_{a_1,b_1}),\s_i(e_{a_2,b_2}),\do,\s_i(e_{a_s,b_s}))$$
if $i$ is even.

\subhead 2.8\endsubhead
Let $\e\in\cz_D$, $\e\ne\emp$.
Let $e_{a,b}\in\e$ be such that $b-a$ is minimum.
If $b-a=0$ we set $i=a=b$; we have $i\in[1,D]$ and $i$ is odd.
If $b-a>0$ we define $i\in[1,D]$ by $a=i-1<i+1\le b$; then $i$ is even..
We will show that

(a) $\e$ is in the image of $\Sigma_i:\cz_{D-2}@>>>\cz_D$.
\nl
If $i$ is odd we can write

$\e=(e_{\ta_1,\tb_1},e_{\ta_2,\tb_2},\do,e_{\ta_s,\tb_s},e_i)$.
\nl
If $i$ is even we can write

$\e=(e_{\ta_1,\tb_1},e_{\ta_2,\tb_2},\do,e_{\ta_s,\tb_s})$
\nl
where $a_t=a,b_t=b$ for some $t$.

To $e_{\ta_t,\tb_t}$ ($t=1,2,\do,s$) we associate the element

$e_{a_t,b_t}=e_{\ta_t-2,\tb_t-2}$ if $i\le\ta_t-2$,

  $e_{a_t,b_t}=e_{\ta_t,\tb_t-2}$ if  $\ta_t<i\le\tb_t-1$,

$e_{a_t,b_t}=e_{\ta_t,\tb_t}$ if $\tb_t<i$.

(Note that we cannot have $i=\ta_t$ or $i=\tb_t$. Moreover
when $i$ is even we see from the definitions that we cannot have
$i=\ta_t-1$.)
This element is in $\cz^*_{D-2}$.

Consider $n\ne m$ in $\{1,2,\do,s\}$. We set

$(\ta_n,\tb_n,\ta_m,\tb_m)=(\ta,\tb,\ta',\tb')$

$(a_n,b_n,a_m,b_m)=(a,b,a',b')$.

We show:

(i) If $\tb<\ta'$, then  $b<a'$.

(ii) If $\ta'<\ta\le\tb<\tb'$, then $a'<a\le b<b'$.
\nl
In the setup of (i) assume that $a'\le b$. Then we have 
$\ta'\le\tb$ or $\ta'-2\le\tb$ or $\ta'-2\le\tb-2$ or $\ta'\le\tb-2$.
The first 3 cases are clearly impossible. In the 4th case we have
$\tb<\ta'\le\tb-2$ hence $\tb<\tb-2$ a contradiction. Thus $b<a'$.

In the setup of (ii), $a',a,b,b'$ is as follows:

$\ta'-2,\ta-2,\tb-2,\tb'-2$  if $i\le\ta'-2$;

$\ta',\ta-2,\tb-2,\tb'-2$  if $\ta'<i\le\ta-2$ (so that $\ta'<\ta-2$);

$\ta',\ta,\tb-2,\tb'-2$  if $\ta<i\le\tb-1$ (so that $\ta\le\tb-2$);
                           
$\ta',\ta,\tb,\tb'-2$  if $\tb<i\le\tb'-2$ (so that $\tb<\tb'-2$);
                            
$\ta',\ta,\tb,\tb'$ if $\tb'<i$.
\nl
Since $i$ is distinct from each of
$\ta',\ta'-1,\ta,\ta-1,\tb,\tb',\tb'-1$ we see that
we must be in one of the 5 cases above. Note that
$a'<a\le b<b'$ in each case.

From (i),(ii) we see that

$\e':=(e_{a_1,b_1},e_{a_2,b_2},\do,e_{a_s,b_s})$ belongs to
$\cz_{D-2}$.
\nl
From the definitions we see that $\e=\Sigma_i(\e')$. Hence (a) holds.

\subhead 2.9\endsubhead
We define a subset $\cz'_D$ of $\cz_D$ by induction on $D$. If $D=0$,
$\cz'_D$ consists of the empty sequence.
Assume now that $D\ge2$. A sequence  $\e\in\cz_D$ is said to be in
$\cz'_D$ if either $\e$ is the empty sequence or if there exists
$i\in[1,D]$ and $\e'\in\cz'_{D-2}$ such that $\e=\Sigma_i(\e')$.
(Note that $\Sigma_i(\e')$ is well defined.) Using 2.8(a)
we see by induction on $D$ that

(a) $\cz_D=\cz'_D$.

\subhead 2.10\endsubhead
Assume that $D\ge2$ and $i\in[1,D]$. For $\e'\in\cz_{D-2}$ we have

(a) $\Th_D(\Sigma_i(\e'))=\ct_i(\Th_{D-2}\e')+Fe_i$ if $i$ is odd;

(b) $\Th_D(\Sigma_i(\e'))=\ct_i(\Th_{D-2}\e')$ if $i$ is even.

We can write $\e'=(e_{a_1,b_1},e_{a_2,b_2},\do,e_{a_s,b_s})$. Then
$$\Th_D(\Sigma_i(\e'))=
F\s_i(e_{a_1,b_1})+F\s_i(e_{a_2,b_2})+\do+F\s_i(e_{a_s,b_s})+Fe_i$$
if $i$ is odd,
$$\Th_D(\Sigma_i(\e'))=
F\s_i(e_{a_1,b_1})+F\s_i(e_{a_2,b_2})+\do+F\s_i(e_{a_s,b_s})$$
if $i$ is even.

Using the definitions we see that
$$\align&\Th_D(\Sigma_i(\e'))=
F\ct_i(e_{a_1,b_1})+F\ct_i(e_{a_2,b_2})+\do+F\ct_i(e_{a_s,b_s})+Fe_i\\&
=\ct_i(Fe_{a_1,b_1}+Fe_{a_2,b_2}+\do+Fe_{a_s,b_s})+Fe_i=
\ct_i(X_{D-1}(\e'))+Fe_i\endalign$$
if $i$ is odd,
$$\align&\Th_D(\Sigma_i(\e'))=
F\ct_i(e_{a_1,b_1})+F\ct_i(e_{a_2,b_2})+\do+F\ct_i(e_{a_s,b_s})\\&
=\ct_i(Fe_{a_1,b_1}+Fe_{a_2,b_2}+\do+Fe_{a_s,b_s})=\ct_i(X_{D-1}(\e'))
\endalign$$
if $i$ is even.
This proves (a),(b).

\subhead 2.11\endsubhead
We prove the following part of Theorem 2.4.

(a) The map $\Th_D$ in 2.4 is well defined.
\nl
We argue by induction on $D$. When $D=0$, (a) is obvious. Assume now
that $D\ge2$.
Let $\e\in\cz_D$. If $\e=\emp$ then $\Th_D(\e)=0\in\cf_D$.
Assume now that $\e\ne\emp$. Using 2.8, we can find $i\in[1,D]$
and $\e'\in\cz_{D-2}$ such that $\e=\Sigma_i(\e')$ so that
$\Th_D(\e)=\Th_D(\Sigma_i(\e'))$.
By the induction hypothesis we have
$\Th_{D-2}\e'\in\cc(V^1_{D-2})$.
By the definition of $\cc(V^1_D)$ we then have

$\ct_i(\Th_{D-2}\e')+Fe_i\in\cc(V^1_D)$ if $i$ is odd;
$\ct_i(\Th_{D-2}\e')\in\cc(V^1_D)$ if $i$ is even.
\nl
Using 2.10, we can rewrite this as $\Th_D(\e)\in\cc(V^1_D)$.
This proves (a).

\subhead 2.12\endsubhead
We prove the following part of Theorem 2.4.

(a) The map $\Th_D$ in 2.4 (see 2.11(a)) is surjective.
\nl
We argue by induction on $D$. When $D=0$, (a) is obvious. Assume now
that $D\ge2$. Let $\ce\in\cc(V^1_D)$. If $\ce=0$ then
$\ce=\Th_D(\emp)$.
Assume now that $\ce\ne0$. We can find $i\in[1,D]$ and
$\ce'\in\cc(V^1_{D-2})$
such that $\ce=\ct_i(\ce')+Fe_i$ if $i$ is odd and
$\ce=\ct_i(\ce')$ if $i$ is even.
By the induction hypothesis we have $\ce'=\Th_{D-2}(\e')$ for some
$\e'\in\cz_{D-2}$. Thus we have $\ce=\ct_i(\Th_{D-2}\e')+Fe_i$ if
$i$ is odd, $\ce=\ct_i(\Th_{D-2}\e')$ if $i$ is even.
Using 2.10 we can rewrite this as $\ce=\Th_D(\e)$ where
$\e=\Sigma_i(\e')\in\cz_D$. This proves (a).

\subhead 2.13\endsubhead
We have $\cc(V^1_D)=\sqc_{s\in[0,d]}\cc^s(V^1_D)$
where $\cc^s(V^1_D)=\{\ce\in\cc(V^1_D);\dim\ce=s\}$.
Clearly, the map $\Th$ in 2.4 restricts for any $s\in[0,d]$ to
a map $\Th^s:\cz^s_D@>>>\cc^s(V^1_D)$.
From 2.12(a) it follows that $\Th^s$ is surjective for any $s\in[0,d]$.
In \cite{HM} it is shown that 
$|\cz^s_D|=N_{d+1,s+1}$ (see 0.5) for any $s\in[0,d]$. Using this and
0.5(a) we see that
$$Cat_{d+1}=\sum_{s\in[0,d]}N(d+1,s+1)=\sum_{s\in[0,d]}|\cz^s_D|=
|\cz_D|.$$
We see that $\Th_D$ is a surjective map from a set with
cardinal $|\cz_D|=Cat_{d+1}$ to a set with the same cardinal
$|\cc(V^1_D)|=|\cf^0_{D/2}(V_D)|=Cat_{d+1}$ (the first equality
holds by 2.2(a); the second equality follows from Theorem 0.4).
It follows that $\Th$ is a bijection and Theorem 2.4 is proved.

This implies that $\Th^s:\cz^s_D@>>>\cc^s(V^1_D)$ is a bijection for
any $s\in[0,d]$. We see that Theorem 0.6 holds. (We use that $X_c$ in
0.5 is the same as $\cc^s(V^1_D)$ if we identify$V^1_D=F^d$.)

\head 3. A conjecture\endhead
\subhead 3.1\endsubhead
In this section we fix a unipotent element $u\in G$. We assume that
either

$G$ is of type $C_{d(d+1)},d\ge1$ and $u$ has
Jordan blocks of sizes $2d,2d,2d-2,2d-2,\do,2,2$ or that

$G$ is of type $B_{d(d+1)},d\ge1$ and $u$ has
Jordan blocks of sizes $2d+1,2d-1,2d-1,\do,1,1$. 

Let $\cb_u$ be the variety of Borel subgroups of $G$ that contain $u$
and let $[\cb_u]$ be the set of irreducible components of $\cb_u$.
Let $A(u)$ be the group of components of the centralizer of $u$ in $G$.
Note that $A(u)$ acts naturally by permutations on $[\cb_u]$.
For each $\x\in[\cb_u]$ we denote by $A(u)_\x$ the stabilizer of $\x$
in $A(u)$. Let $\Xi_u$ be the set of subgroups of $A(u)$ of the form
$A(u)_\x$ for some $\x\in[\cb_u]$.

We assume that $c$ is the family containing the Springer representation
of $W$ associated to $u$ and to the unit representation of $A(u)$.
We conjecture that

(a) {\it there exists an isomorphism $A(u)@>\si>>V^1_D,D=2d$ which
carries $\Xi_u$ to the collection $\cc(V^1_D)$ (see 2.1) of subspaces
of $V^1_D$.}
\nl
(This would imply that $|\Xi_u|$ is a Catalan number.)

We have verified that (a) is true when $d=1,2,3$.

\enddocument